\documentstyle[amstex,amssymb,12pt,leqno]{article}
\newcommand{\df}{\sc}

\newcommand{\R}{\mbox{$I\hspace{-0.3em}R$}}
\newcommand{\Proj}{\mbox{$I\hspace{-0.3em}P$}}

\newcommand{\g}{\frak{g}}%{\mbox{\boldmath$g$}}
\newcommand{\p}{\frak{p}}%{\mbox{\boldmath$p$}}
\newcommand{\k}{\frak{k}}%{\mbox{\boldmath$k$}}
\newcommand{\ol}{\frak{o}}%{\mbox{\boldmath$o$}}
\renewcommand{\a}{\frak{a}}
\newcommand{\tr}{\mbox{tr}}
\newcommand{\Ad}{\mbox{Ad}}

\newcommand{\bM}[1]{
\renewcommand{\baselinestretch}{1}\footnotesize$\left(\begin{array}
{*{#1}{c}} }
\newcommand{\eM}{ \end{array}\right)$\renewcommand{\baselinestretch}{1}
\normalsize }
\newcommand{\bMM}[1]{ \renewcommand{\baselinestretch}{1}\scriptsize$
\left(\begin{array}{*{#1}{c}} }
\newcommand{\eMM}{
\end{array}\right)$\renewcommand{\baselinestretch}{1}\normalsize }
\newcommand{\bBM}{
\renewcommand{\baselinestretch}{1}\small$\left(\begin{array}{cc} }
\newcommand{\eBM}{
\end{array}\right)$\renewcommand{\baselinestretch}{1}\normalsize }
\newcommand{\bBl}[1]{
\renewcommand{\baselinestretch}{1}\scriptsize$\begin{array}
{*{#1}{c}} }
\newcommand{\eBl}{ \end{array}$\renewcommand{\baselinestretch}{1}\normalsize }
\newcommand{\bdr}{ \renewcommand{\baselinestretch}{1.4}\scriptsize$
\left\{\begin{array}{c} }
\newcommand{\edr}{
\end{array}\right.$\renewcommand{\baselinestretch}{1}\normalsize }

\author{F.~Burstall, U.~Hertrich-Jeromin, F.~Pedit\thanks{Partially supported
by NSF grant DMS 2905293}, U.~Pinkall}
\title{{\bf Curved Flats and Isothermic Surfaces}}
\begin{document}
\maketitle
\begin{abstract}
We show how pairs of isothermic surfaces are given by curved flats in a
pseudo Riemannian symmetric space and vice versa. Calapso's fourth order
partial differential equation is derived and, using a solution of this
equation, a M\"obius invariant frame for an isothermic surface is built.
\end{abstract}

\section{Introduction}\label{1}
These notes grew out of a series of discussions on a recent paper by
J.~Cie\'{s}li\'{n}ski, P.~Goldstein and A.~Sym \cite{CGS}: these
authors give a characterization of isothermic surfaces as "soliton
surfaces" by introducing a spectral parameter. In trying to
understand the geometric meaning of this spectral parameter, we
observed some analogies with the theory of conformally flat hypersurfaces
in a four-dimensional space form: Guichard's nets may be
understood as a kind of analogue of isothermic parametrizations of
Riemannian surfaces (cf.\cite[no.3.4.1]{J}), and so it seems
natural to look for relations between the theory of isothermic
surfaces in three-dimensional space forms and the theory of
conformally flat hypersurfaces in four-dimensional space forms.  Here
we would like to present some results we found --- especially the
possibility of constructing isothermic surfaces using
\section{Curved Flats}\label{2}
A curved flat is the natural generalization of a developable surface
in Euclidean space: it is a submanifold $M\subset G/K$ of a
(pseudo-Riemannian) symmetric space for which the curvature operator
of $G/K$ vanishes on $\bigwedge^2 TM$\footnote{Thus $M$ is {\em
curvature isotropic\/} in the sense of \cite{FP}}.  Thus, a curved
flat may be thought of as the enveloping submanifold of a congruence
of flats --- totally geodesic submanifolds --- of the symmetric space.
Taking a regular parametrization $\gamma:M\rightarrow G/K$ of a curved
flat and a framing $F:M\rightarrow G$ of this parametrization, the
Maurer-Cartan form $\Phi=F^{-1}dF$ of the framing has a natural
decomposition $\Phi=\Phi_{\k}+\Phi_{\p}$ according to the symmetric
decomposition\footnote{Thus $\k$ and $\p$ are the $+1$ and
$-1$-eigenspaces, respectively, of the involution fixing $\k$ and so
satisfy the characteristic conditions \[ [\k,\k]\subset\k,
[\k,\p]\subset\p, [\p,\p]\subset\k. \]} $\g=\k\oplus\p$ of the Lie
algebra $\g$. Now the condition for $\gamma$ to parametrize a curved
flat may be formulated as\footnote{The product
\[ [\Phi\wedge\Psi](v,w):=[\Phi(v),\Psi(w)]-[\Phi(w),\Psi(v)] \]
defines a symmetric product on the space of Lie algebra valued 1-forms
with values in the space of Lie algebra valued 2-forms.}
\begin{equation}\label{200}
   [ [\Phi_{\p}\wedge\Phi_{\p}],\p]\equiv 0\ .
\end{equation}
In case that $G$ is semisimple, it is straightforward\footnote{In
fact, $[\p,\p]\oplus\p$ is an ideal of $\g$ so that we have a
decomposition $\g=\k'\oplus[\p,\p]\oplus\p$ where $\k'$ is a
complementary ideal commuting with $[\p,\p]\oplus\p$.  Thus, if
$\a\subset\p$ satisfies $[[\a,\a],\p]=0$ we deduce that $[\a,\a]$ lies
in the center of $\g$ and so vanishes.} to see that this is
equivalent to
\begin{equation}\label{201}
    [\Phi_{\p}\wedge\Phi_{\p}]\equiv 0\ .
\end{equation}

To summarise, we have the
\\ {\bf Definition} of a curved flat: An immersion
$\gamma:M\rightarrow G/K$ is said to parametrize a {\df curved
flat}, if the $\p$-part in the symmetric decomposition
 of the Maurer-Cartan form
$F^{-1}dF= \Phi=\Phi_{\k}+\Phi_{\p}$ of a framing $F:M\rightarrow G$ of
$\gamma$ defines a congruence $p\mapsto\Phi_{\p}|_p(T_pM)$ of
abelian subalgebras of $\g$.

At this point we should remark that curved flats naturally arise in one
parameter families \cite{FP}: setting
\begin{equation}\label{202}
    \Phi_{\lambda}:=\Phi_{\k}+\lambda\Phi_{\p}
\end{equation}
the Maurer-Cartan equation
$d\Phi_{\lambda}+\frac{1}{2}[\Phi_{\lambda}\wedge\Phi_{\lambda}]=0$
for the loop $\lambda\mapsto\Phi_{\lambda}$ of forms splits into the three
equations
\begin{equation}\label{203}
\begin{array}{l}
    0 = d\Phi_{\k} + \mbox{$\frac{1}{2}$}[\Phi_{\k}\wedge\Phi_{\k}] \cr
    0 = d\Phi_{\p} + [\Phi_{\k}\wedge\Phi_{\p}] \cr
    0 = [\Phi_{\p}\wedge\Phi_{\p}]\ , \cr
\end{array}
\end{equation}
and hence the integrability of the loop $\lambda\mapsto\Phi_{\lambda}$ is
equivalent to the forms $\Phi_{\lambda}$ being the Maurer-Cartan forms for
some framings $F_{\lambda}:M\rightarrow G$ of curved flats
$\gamma_{\lambda}:M\rightarrow G/K$.  Thus integrable systems theory
may be applied to produce examples.

Now we will consider the case leading to the theory of isothermic surfaces:
let
\begin{equation}\label{204}
    G := O_1(5)\hspace{1em}\mbox{and}\hspace{1em}
    K := O(3)\times O_1(2)\ .
\end{equation}
The coset space $G_+(5,3)=G/K$ of space-like 3-planes in the Minkowski
space $\R^5_1$ becomes a six dimensional pseudo-Riemannian symmetric
space of signature $(3,3)$ when endowed with the metric induced by the
Killing form.  We will
consider two-dimensional curved flats
\begin{equation}\label{205}
    \gamma:M^2\rightarrow G_+(5,3)\
\end{equation}
satisfying the regularity assumption that the metric on $M^2$ induced
by $\gamma$ is non-degenerate.

Fixing a pseudo orthonormal basis $(e_1,\dots,e_5)$ of the Minkowski
space $\R^5_1$ with
\begin{equation}\label{206}
    (\langle e_i,e_j\rangle)_{ij} = E_5 :=
    \mbox{\bBM I_3 & 0 \cr 0 & \mbox{\bBl{2} 0&1\cr 1&0\cr \eBl}\cr
    \eBM} \ ,
\end{equation}
we get the matrix representations
\begin{equation}\label{207}
\begin{array}{l}
    O_1(5) = \{A\in Gl(5,\R)| A^tE_5A = E_5\} \cr
    \ol_1(5) = \{X\in \frak{gl}(5,\R)| (E_5X)+(E_5X)^t=0\}\ .\cr
\end{array}
\end{equation}
The subalgebra $\k$ and its complementary linear subspace $\p$ in the
symmetric decomposition of $\ol_1(5)$ are given by the $+1$-
resp. $-1$-eigenspaces of the involutive automorphism
$\Ad(Q):\ol_1(5)\rightarrow\ol_1(5)$ with $Q=\mbox{\bMM{2} -I_3&0\cr
0&I_2\cr\eMM}$. Writing down the Maurer-Cartan form of a framing
$F:M^2\rightarrow O_1(5)$ of our curved flat $\gamma:M^2\rightarrow
G_+(5,3)$ with this notation we obtain
\begin{equation}\label{208}
\begin{array}{ll}
    F^{-1}dF = \Phi = \Phi_{\k}+\Phi_{\p} & \mbox{with} \cr
    \Phi_{\k} = \mbox{\bM{2} \Omega & 0 \cr 0 & \mbox{\large$\nu$} \cr\eM}: &
        TM\rightarrow\ol(3)\times\ol_1(2) \cr
    \Phi_{\p} = \mbox{\bM{2} 0 & \mbox{\large$\eta$} \cr
	-E_2\mbox{\large$\eta$}^t
        & 0 \cr\eM}: & TM\rightarrow\p\ . \cr
\end{array}
\end{equation}
The image of $\Phi_{\p}$ at each $p\in M^2$ is a $2$-dimensional
abelian subspace of $\p$ on which the Killing form is non-degenerate.
One can show that there are precisely two $K$-orbits of maximal
abelian subspaces of $\p$: one consists of $3$-dimensional
subspaces which are isotropic for the Killing form while the other
consists of $2$-dimensional subspaces on which the Killing form has
signature $(1,1)$.  We therefore conclude that the images of each
$\Phi_{\p}$  are maximal abelian and $K$-conjugate and so we can put
\mbox{\large$\eta$} into the standard form
\begin{equation}\label{209}
    \mbox{\Large$\eta$} = \mbox{\bMM{2} \omega_1 & -\omega_1 \cr
        \omega_2 & \omega_2 \cr 0 & 0 \cr\eMM}\hspace{1em}
\end{equation}
by applying a gauge transformation $M\to K$.

Calculating the Maurer-Cartan equation using the ansatz
\begin{equation}\label{210}
    \Omega = \mbox{\bMM{3} 0 & \omega & -\psi_1 \cr -\omega & 0 & -\psi_2 \cr
        \psi_1 & \psi_2 & 0 \cr\eMM}\hspace{1em}\mbox{and}\hspace{1em}
    \mbox{\Large$\nu$} = \mbox{\bM{2} \nu & 0 \cr 0 & -\nu \cr\eM}
\end{equation}
together with \mbox{\large$\eta$} given by (\ref{209}), we see that
\begin{equation}\label{211}
    d\omega_1 = d\omega_2 = 0\ .
\end{equation}
So we are given canonical coordinates $(x,y):M\rightarrow\R^2$ by integrating%
\footnote{Since our theory is local, all closed forms may be assumed to be
exact.} the forms $\omega_1$ and $\omega_2$. Moreover, since
we also have $d\nu=0$, we may set $\nu=-du$ for a suitable function $u\in
C^{\infty}(M)$ --- this gives us $\omega=u_ydx-u_xdy$, where $u_x$ and $u_y$
denote the partial derivatives of $u$ in $x$- resp. $y$-directions.
Finally, the equations $\psi_1\wedge\omega_1=0$ and $\psi_2\wedge\omega_2=0$
show that $\psi_1=e^uk_1dx$ and $\psi_2=e^uk_2dy$ for two functions
$k_i\in C^{\infty}(M)$.

We now perform a final $O_1(2)$-gauge $\mbox{\bBM I_3&0\cr
0&\mbox{\bBl{2} e^u&0\cr 0&e^{-u}\cr\eBl}\cr\eBM}:M\rightarrow
O(3)\times O_1(2)$ and insert the spectral parameter $\lambda$ to
obtain the Maurer-Cartan form discussed in (cf.\cite{CGS}):
\begin{equation}\label{212}
    \Phi_{\lambda} = \mbox{\bMM{6} 0 & u_ydx-u_xdy & -e^uk_1dx && \lambda
	e^u dx& -\lambda e^{-u}dx \cr
        -u_ydx+u_xdy & 0 & -e^uk_2dy && \lambda e^udy & \lambda e^{-u}dy \cr
        e^uk_1dx & e^uk_2dy & 0 && 0 & 0 \vspace{1ex}\cr
        \lambda e^{-u}dx & -\lambda e^{-u}dy & 0 && 0 & 0 \cr
        -\lambda e^udx & -\lambda e^udy & 0 && 0 & 0 \cr\eMM}\ .
\end{equation}
We are now lead directly to the theory of
\section{Isothermic Surfaces}\label{3}
In the context of M\"obius geometry the three sphere $S^{3}$ is
viewed as the projective light-cone $\Proj L^{4}$ in $\R^{5}_{1}$ while
the Lorentzian sphere $\{v\in\R^5_1| \langle v,v\rangle=1\}$ should
be interpreted as the space of (oriented) spheres in the three
sphere\footnote{Or, equivalently, it may be interpreted as the space
of (oriented) spheres and planes in Euclidean three space $\R^3$: the
polar hyperplane to a vector $v$ of the Lorentz sphere intersects the
three sphere --- thought of as the absolute quadric in projective four
space --- in a two sphere. Stereographic projection yields a sphere
in $\R^3$ or, if the projection center lies on the sphere, a plane.}
(cf.\cite{B}).
Now, denoting by
\begin{equation}\label{301}
\begin{array}{l}
    n := Fe_3:M\rightarrow S_1^5=\{v\in\R_1^5|\langle v,v\rangle=1\} \cr
    f := Fe_4:M\rightarrow L^4=\{v\in\R_1^5|\langle v,v\rangle=0\} \cr
    \hat{f} := Fe_5:M\rightarrow L^4 \cr
\end{array}
\end{equation}
one of the sphere congruences resp. the two immersions given by our frame
$F$, we see that \\
{\bf Theorem:} {\em The sphere congruence $n$ given by our curved flat is a
Ribeaucour sphere congruence\footnote{\rm The curvature lines on the two
enveloping immersions correspond.}, which is enveloped by two isothermic
immersions $f$ and $\hat{f}$} (cf.\cite[p.362]{B}):\\
Since
\begin{equation}\label{302}
\begin{array}{l}
    \langle f,n\rangle = 0\hspace{1em}\mbox{and}\hspace{1em}
        \langle df,n\rangle \equiv 0,\cr
    \langle \hat{f},n\rangle = 0\hspace{1em}\mbox{and}\hspace{1em}
        \langle d\hat{f},n\rangle \equiv 0\ ,\cr
\end{array}
\end{equation}
the immersions $f$ and $\hat{f}$ do envelop the sphere congruence $n$ and,
since the bilinear forms
\begin{equation}\label{303}
\begin{array}{l}
    \langle df,dn\rangle = \lambda e^{2u}(k_1dx^2+k_2dy^2),\cr
    \langle d\hat{f},dn\rangle = \lambda (-k_1dx^2+k_2dy^2)\cr
\end{array}
\end{equation}
are diagonal with respect to the induced metrics
\begin{equation}\label{304}
\begin{array}{l}
    \langle df,df\rangle = \lambda^2e^{2u}(dx^2+dy^2),\cr
    \langle d\hat{f},d\hat{f}\rangle = \lambda^2e^{-2u}(dx^2+dy^2)\ ,\cr
\end{array}
\end{equation}
the two immersions $f$ and $\hat{f}$ are isothermic\footnote{The bundle
defined by $\mbox{span}\{n,f,\hat{f}\}$ over $M$ is flat (cf.(\ref{212})) and
so
the map $p\mapsto d_pf(T_pM)$ defines a normal congruence of
circles \cite{Cool}: for each $p\in M$ \[t\mapsto f_t(p):=
\frac{1}{\sqrt{2}}\sin t\cdot n(p)+\frac{1}{2}(1+\cos t)\cdot
f(p)-\frac{1}{2}(1-\cos t)\cdot\hat{f}(p)\]
parametrizes the circle $(d_pf(T_pM))^{\perp}$ meeting the sphere $n(p)$ in
$f(p)$ and $\hat{f}(p)$ orthogonal. Since $n$, $f$ and $\hat{f}$ are parallel
sections in this bundle, the maps $p\mapsto f_t(p)$ (which generically are
not degenerate) parametrize the surfaces orthogonal to this congruence of
circles.

In general the immersions $f$ and $\hat{f}=f_{\pi}$ will be the only
isothermic surfaces among the surfaces.\label{norm}}.

It is quite difficult to calculate the first and second fundamental forms
of these isothermic immersions, when they are projected to $S^3$ resp. $\R^3$,
but applying a (constant) conformal change (constant $O_1(2)$-gauge)
\begin{equation}\label{305}
\begin{array}{l}
    f\leadsto\mbox{$\frac{1}{\lambda}$}f\hspace{1em}\mbox{and}\hspace{1em}
        \hat{f}\leadsto\lambda\hat{f}\hspace{1em}\mbox{or} \cr
    f\leadsto\lambda f\hspace{1em}\mbox{and}\hspace{1em}
        \hat{f}\leadsto\mbox{$\frac{1}{\lambda}$}\hat{f} \cr
\end{array}
\end{equation}
and sending $\lambda\rightarrow 0$, $\hat{f}$ resp. $f$ becomes a
constant vector --- $\Phi_{\lambda=0}e_5$ resp. $\Phi_{\lambda=0}e_4$
vanishes. This constant light-like vector may be interpreted as the
point at infinity and we therefore obtain an isothermic immersion
$f:M\rightarrow\R^3$ with first and second fundamental forms
\begin{equation}\label{306}
\begin{array}{l}
    I = e^{2u}(dx^2+dy^2) \cr
    II = e^{2u}(k_1dx^2+k_2dy^2) \cr
\end{array}
\end{equation}
resp. its Euclidean dual surface $\hat{f}:M\rightarrow\R^3$ with first and
second fundamental forms
\begin{equation}\label{307}
\begin{array}{l}
    \hat{I} = e^{-2u}(dx^2+dy^2) \cr
    \hat{II} = -k_1dx^2+k_2dy^2\ . \cr
\end{array}
\end{equation}
We now recognise the remaining three equations from the Maurer-Cartan
equation for $\Phi_{\lambda}$
\begin{equation}\label{308}
\begin{array}{l}
    0 = \Delta u +e^{2u}k_1k_2 \cr
    0 = k_{1y} + (k_1-k_2)u_y \cr
    0 = k_{2x} - (k_1-k_2)u_x \cr
\end{array}
\end{equation}
as the Gau\ss{} and Codazzi equations of the Euclidean immersion $f$
resp.  its dual $\hat{f}$ \footnote{The Euclidean dual of an isothermic
surface is obtained by integrating the closed 1-form
$d\hat{f}:=e^{-2u}(-f_xdx+f_ydy)$: see for example \cite[p.14]{S}.

When the normal congruence of
circles mentioned in footnote \ref{norm} %(p.\pageref{norm})
is projected to Euclidean three space $\R^3$, we see that, in the limit
$\lambda\rightarrow 0$, the circles become straight lines --- circles
meeting the collapsed surface $\hat{f}$ resp. $f$ in the point at
infinity --- while the Ribeaucour sphere congruence enveloped by the
two surfaces $f$ and $\hat{f}$ becomes the congruence of tangent
planes of $f$ resp. $\hat{f}$.} \cite{C}, \cite{CGS}.  As a
consequence, we can invert our construction and build a curved flat
from an isothermic surface:\\ {\bf Theorem}. {\em Given an isothermic
surface $f:M^2\rightarrow\R^3$ and its Euclidean dual surface
$\hat{f}:M\rightarrow\R^3$ we get a curved flat $\gamma:M\rightarrow
G_+(5,3)$ integrating the Maurer-Cartan form (\ref{212}), which we
are able to write down knowing the first and second fundamental forms
of the immersions $f$ and $\hat{f}$ \footnote{Since this construction
depends on the conformal rather than the Euclidean geometry of the
ambient space, we generally get a whole three parameter family of
loops of curved flats from one isothermic surface: when viewing
our given isothermic surface as a surface in the three sphere $S^3$,
we may choose the point at infinity arbitrarily.}.}

Another way to obtain these two Euclidean immersions is presented in
\cite{CGS}. Applying Sym's formula to the associated family of frames
$F=F(\lambda)$, we obtain a map
\begin{equation}\label{309}
    (\mbox{$\frac{\partial}{\partial\lambda}$}F)F^{-1}|_{\lambda=0}:
        M\rightarrow\p\ ;
\end{equation}
interpreting $\p$ as two copies of Euclidean three space\footnote{Here the
Euclidean metric is induced by the quadratic form
$\frac{1}{2}\tr\Phi_{\p}^t\Phi_{\p}$
instead of the Killing form.} $\R^3$ this map gives us the immersion $f$, and
in the other copy of $\R^3$, its dual $\hat{f}$: this can be seen by looking
at the differential
\begin{equation}\label{310}
\begin{array}{l}
    d(\mbox{$\frac{\partial}{\partial\lambda}$}F)F^{-1}|_{\lambda=0}
        = F_0\Phi_{\p}F_0^{-1} \vspace{0.5ex}\cr\hspace{1em}
        \cong H_3\mbox{\bMM{2} e^udx & -e^{-u}dx \cr e^udy & e^{-u}dy \cr 0 & 0
\cr\eMM}\ . \cr
\end{array}
\end{equation}
Here $F_0=\mbox{\bM{2} H_3&0\cr 0&I_2\cr\eM}$ solves the equation
$F_0^{-1}dF_0=\Phi_{\k}$ and thus we may view
$H_3:M\rightarrow O(3)$ as a Euclidean framing of $f$ resp. $\hat{f}$.

There is another possibility for producing isothermal surfaces in Euclidean
space $\R^3$ (or $S^3$): that is, by using a solution of
\section{Calapso's equation}\label{4}
To understand this, we write down the Maurer-Cartan form of a frame
$F:M\rightarrow O_1(5)$, which is M\"obius-invariantly connected to a given
immersion: taking $f=Fe_4$ the (unique) isometric lift of the isothermic
immersion and $n=Fe_3$ the central sphere congruence (conformal Gau\ss{}
map) of the immersion, the frame is determined by the assumption of being an
adapted frame (i.e. $Fe_1=f_x$ and $Fe_2=f_y$). The associated Maurer-Cartan
form will be
\begin{equation}\label{401}
    \Phi = \mbox{\bMM{6} 0 & 0 & kdx && dx & \chi_1 \cr
        0 & 0 & -kdy && dy & \chi_2 \cr
        -kdx & kdy & 0 && 0 & \tau \vspace{1ex}\cr
        -\chi_1 & -\chi_2 & -\tau && 0 & 0 \cr
        -dx & -dy & 0 && 0 & 0 \cr\eMM}\ ,
\end{equation}
$k^2$ being the conformal factor relating the metric induced by the central
sphere congruence to the isometric one, and the 1-forms $\chi_1$, $\chi_2$ and
$\tau$ to be determined. From the Maurer-Cartan equation for this form we
learn that
\begin{equation}\label{402}
\begin{array}{l}
    \tau = k_xdx - k_ydy \cr
    \chi_1 = (\mbox{$\frac{1}{2}$}k^2-u)dx - \mbox{$\frac{k_{xy}}{k}$}dy \cr
    \chi_2 = -\mbox{$\frac{k_{xy}}{k}$}dx + (\mbox{$\frac{1}{2}$}k^2+u)dy\ ,
\cr
\end{array}
\end{equation}
where $u\in C^{\infty}(M)$ is a function satisfying the differential equation
\begin{equation}\label{403}
    du = -((\mbox{$\frac{k_{xy}}{k}$})_y+(k^2)_x)dx
        + ((\mbox{$\frac{k_{xy}}{k}$})_x+(k^2)_y)dy
\end{equation}
--- the integrability condition of this equation is a fourth order partial
differential equation closely related to Calapso's original equation
\cite{C}:
\begin{equation}\label{404}
    0 = \Delta(\mbox{$\frac{k_{xy}}{k}$})+2(k^2)_{xy}
\end{equation}
This shows, that\\
{\bf Theorem:} {\em Any isothermic surface gives rise to a solution
of Calapso's equation.

Conversely, from a solution $k\in C^{\infty}(M)$ of Calapso's
equation we can construct a M\"obius invariant frame of an isothermic surface
by integrating the Maurer-Cartan form (\ref{401}), where the function $u$
is a solution of (\ref{403}).}

Now, applying a conformal change $f\leadsto\frac{1}{k}f$ while fixing the
central sphere congruence $n\leadsto n$, the Maurer-Cartan form of
the associated frame becomes
\begin{equation}\label{406}
    \Phi = \mbox{\bMM{6} 0 & \omega & kdx && \frac{1}{k}dx & \chi_1 \cr
        -\omega & 0 & -kdy && \frac{1}{k}dy & \chi_2 \cr
        -kdx & kdy & 0 && 0 & 0 \vspace{1ex}\cr
        -\chi_1 & -\chi_2 & 0 && 0 & 0 \cr
        -\frac{1}{k}dx & -\frac{1}{k}dy & 0 && 0 & 0 \cr\eMM},
\end{equation}
where
\begin{equation}\label{405}
\begin{array}{l}
    \omega = -\mbox{$\frac{k_y}{k}dx+\frac{k_x}{k}dy$} \cr
    \chi_1 = k(\mbox{$\frac{k_{xx}}{k}-
	\frac{k_x^2+k_y^2}{2k^2}+\frac{1}{2}k^2-u$})dx \cr
    \chi_2 = k(\mbox{$\frac{k_{yy}}{k}-
	\frac{k_x^2+k_y^2}{2k^2}+\frac{1}{2}k^2+u$})dy \cr
\end{array}\ .
\end{equation}
Here we see that the central sphere congruence of an isothermic
surface is a Ribeaucour sphere congruence, which actually is a
characterisation of isothermic surfaces (cf.\cite[p.374]{B}), and
hence it has flat normal bundle as a codimension two surface in the
Lorentz sphere $S_1^4$.

In general, the second enveloping surface of the central sphere congruence of
an isothermic surface will not be an isothermic surface and it seems to be
difficult to built a curved flat starting with it. But in
a quite simple case this is possible:
\section{Example}\label{5}
Starting with a surface of revolution
\begin{equation}\label{501}
    f(x,y) = (r(x)\cos y, r(x)\sin y, z(x))\ ,
\end{equation}
the functions $r$ and $z$ satisfying the differential equation
\begin{equation}\label{502}
    r^2 = r^{\prime 2} + z^{\prime 2}\ ,
\end{equation}
i.e. the curve $(r,z)$ being parametrized by arc length (thought
of as a curve in the Poincar\'{e} half plane), we may write down the loop of
Maurer-Cartan forms
\begin{equation}\label{503}
    \Phi_{\lambda} = \mbox{\bMM{6} 0 & -\frac{r^{\prime}}{r}dy &
-\frac{r^{\prime}
z^{\prime\prime}-r^{\prime\prime}z^{\prime}}{r^2}dx && \lambda rdx
& -\frac{\lambda}{r}dx \cr
\frac{r^{\prime}}{r}dy & 0 & -\frac{z^{\prime}}{r}dy && \lambda rdy &
\frac{\lambda}{r}dy \cr
\frac{r^{\prime}z^{\prime\prime}-r^{\prime\prime}z^{\prime}}{r^2}dx &
\frac{z^{\prime}}
{r}dy & 0 && 0 & 0 \vspace{1ex}\cr
        \frac{\lambda}{r}dx & -\frac{\lambda}{r}dy & 0 && 0 & 0 \cr
        -\lambda rdx & -\lambda rdy & 0 && 0 & 0 \cr\eMM}\ ,
\end{equation}
which gives us the immersion $f$ and its dual $\hat{f}$ in the limit
$\lambda\rightarrow 0$.

On the other hand, denoting by $H=\frac{1}{2}(\frac{z^{\prime}}{r^2}+
\frac{r^{\prime}z^{\prime\prime}-r^{\prime\prime}z^{\prime}}{r^3})$
the mean curvature of our surface of rotation, the central sphere
congruence of $f$ is $n+Hf$.  The metric it induces has conformal
factor $k^{2}$ (relative to the metric induced by $f$)
given by
\begin{equation}\label{504}
    k = \mbox{$\frac{1}{2r^2}$}(rz^{\prime}-
        r^{\prime}z^{\prime\prime}+r^{\prime\prime}z^{\prime})\ .
\end{equation}
Since $k_y\equiv 0$, this is obviously a solution of Calapso's equation and
a function $u$ solving (\ref{403}) is $u=\lambda^2-k^2$. So the
Maurer-Cartan form (\ref{401}) becomes
\begin{equation}\label{505}
    \Phi_{\lambda} = \mbox{\bMM{6} 0 & 0 & kdx && dx &
		(\frac{3}{2}k^2-\lambda^2)dx \cr
        0 & 0 & -kdy && dy & (-\frac{1}{2}k^2+\lambda^2)dy \cr
        -kdx & kdy & 0 && 0 & k_xdx \vspace{1ex}\cr
        -(\frac{3}{2}k^2-\lambda^2)dx & (\frac{1}{2}k^2-\lambda^2)dy & -k_xdx
&& 0 & 0 \cr
        -dx & -dy & 0 && 0 & 0 \cr\eMM}\ .
\end{equation}
A change $n\leadsto n+kf$ of the sphere congruence, enveloped by $f$, followed
by an $O_1(2)$-gauge $f\leadsto\lambda f$ and
$\hat{f}\leadsto\lambda^{-1}\hat{f}$
gives us the Maurer-Cartan form
\begin{equation}\label{506}
    \Phi_{\lambda} = \mbox{\bMM{6} 0 & 0 & 2kdx && \lambda dx & -\lambda dx \cr
        0 & 0 & 0 && \lambda dy & \lambda dy \cr
        -2kdx & 0 & 0 && 0 & 0 \vspace{1ex}\cr
        \lambda dx & -\lambda dy & 0 && 0 & 0 \cr
        -\lambda dx & -\lambda dy & 0 && 0 & 0 \cr\eMM}\ .
\end{equation}
of a curved flat, quite different from that coming from (\ref{503}).

To understand the geometry of the two enveloping immersions $f=Fe_4$ and
$\hat{f}=Fe_5$, we remark that the sphere congruence $n=Fe_3$ depends only
on one variable and hence the two immersions parametrize a channel surface;
moreover all spheres of the congruence are perpendicular to the fixed circle%
\footnote{We have $\Phi e_2=-(e_4+e_5)dy$ and $\Phi (e_4+e_5)=2e_2dy$.}
given by span$\{Fe_2,F(e_4+e_5)\}$, which may be thought as an axis of
rotation: the immersions $f$ and $\hat{f}$ parametrize pieces of a surface of
revolution\footnote{The meridian curve is given by
$\frac{1}{\sqrt{2}}(f-\hat{f})$
--- which only depends on one variable --- thought as a curve in the
Poincar\'{e}
half plane; its tangent field is given by $Fe_1$ and its unit normal field
by $n=Fe_3$.}, $f$ and $\hat{f}$ being axisymmetric\footnote{The circles
$\{F(p)e_1,F(p)e_2\}^{\perp}$ intersecting the sphere $n(p)$ orthogonally
in $f(p)$ and $\hat{f}(p)$ all meet the axis (cf. footnote \ref{norm}, page
\pageref{norm}).}. Taking now the limit
$\lambda\rightarrow 0$, we obtain a cylinder resp. its dual, which is an
(axial) reflection of the original cylinder.

%\newpage


\vspace{10ex}
\begin{thebibliography}{99}
\bibitem{B}
        {\em W.~Blaschke}, Vorlesungen \"uber Differentialgeometrie III,
        Berlin (1929)
\bibitem{S}
	{\em A.~I.~Bobenko}, Surfaces in Terms of 2 by 2 Matrices. Old
	and New Integrable Cases, SFB~288 Preprint No.~{\bf 66}
\bibitem{C}
        {\em P.~Calapso}, Sulla superficie a linee di curvatura isoterme,
        Rend. Circ. Mat. Palermo~{\bf 17} (1903) 275-286
\bibitem{CGS}
        {\em J.~Cie\'{s}li\'{n}ski, P.~Goldstein, A.~Sym}, Isothermic
        Surfaces in $E^3$ as Soliton Surfaces, Short Report (1994)
\bibitem{Cool}
        {\em J.~Coolidge}, Congruences and Complexes of Circles,
        Trans. AMS~{\bf 15} (1914) 107-134
\bibitem{FP}
        {\em D.~Ferus, F.~Pedit}, Curved Flats in Symmetric Spaces,
        Manuscript (1993)
\bibitem{J}
        {\em U.~Hertrich-Jeromin}, \"Uber konform flache Hyperfl\"achen
        in vierdimensionalen Raumformen, Thesis (1994)
\end{thebibliography}
\end{document}